\RequirePackage{fix-cm}
\documentclass[smallextended]{svjour3}       % onecolumn (second format)
\smartqed  % flush right qed marks, e.g. at end of proof
\usepackage{graphicx}
\usepackage{amsmath,amsfonts,bm}
\usepackage{algorithm}
\usepackage{algpseudocode}
\usepackage{multirow}
\usepackage{xcolor}
\usepackage{tcolorbox}

\begin{document}

\title{A New Continuous Optimization Method for Mixed Integer Space Travelling Salesman Problem %\thanks{Grants or other notes
%about the article that should go on the front page should be
%placed here. General acknowledgments should be placed at the end of the article.}
}
%\subtitle{}

%\titlerunning{Short form of title}        % if too long for running head

\author{Liqiang Hou   \and
        Shufan Wu     \and
        Zhongcheng Mu \and 
        Meilin Liu  %etc.
}

%\authorrunning{Short form of author list} % if too long for running head

\institute{Liqiang Hou \at
              School of Aeronautics and Astronautics, Shanghai Jiaotong University, Shanghai, China, 200240 \\
              \email{houliqiang@sjtu.edu.cn}           %  \\
%             \emph{Present address:} of F. Author  %  if needed
           \and
           Shufan Wu, Zhongcheng Mu, Meilin Liu \at
           School of Aeronautics and Astronautics, Shanghai Jiaotong University, Shanghai, China, 200240
}

\date{Received: date / Accepted: date}
% The correct dates will be entered by the editor

\maketitle

\begin{abstract}
The travelling  salesman problem (TSP) of space trajectory design  is complicated  by its complex structure design space. The graph based tree search and stochastic seeding combinatorial approaches are commonly employed to tackle the time-dependent TSP due to their combinatorial nature. In this paper, a new continuous optimization strategy for the mixed integer  combinatorial problem is proposed. The space trajectory combinatorial problem is tackled using continuous gradient based method. A continuous mapping technique is developed to map the integer type ID of targets on the sequence to a set of continuous design variables. Expected flyby targets are introduced as references and used as priori to select the candidate target to fly by. Bayesian based analysis is employed to model accumulated posterior of the sequence and a new objective function with quadratic form constraints is constructed. The new introduced auxiliary design variables of expected targets together with the original design variables are set to be optimized. A gradient based optimizer is used to search optimal sequence parameter. Performances of the proposed  algorithm are demonstrated through a multiple debris rendezvous problem and a static TSP benchmark.
% \PACS{PACS code1 \and PACS code2 \and more}
% \subclass{MSC code1 \and MSC code2 \and more}
\end{abstract}

\section{Introduction}
\label{intro}
Combinatorial optimization of space trajectory design finds an optimal sequence of the objects from a finite set of objects. The searching process of combinatorial optimization operates on discrete feasible solutions. Typical problems involving combinatorial optimization are the travelling salesman problem (TSP), the minimum spanning tree problem (MST), and the knapsack problem. Different to the common design problems of the space design, in combinatorial optimization of trajectory design a number of sequences are selected and optimized. For each sequence, a global optimization problem is solved in order to search the optimal continuous parameters.

The optimization considers not only the continuous design parameters, e.g. launch time, fuel consumption, impulsive velocity increment, time of flight of the transfers between objects,etc.,  but the discrete parameters,such as the celestial bodies and debris to flyby or rendezvous. One needs to explore not only the continuous design space to search the optimal orbital transfer parameters, but the set of discrete items of the targets to be visited. The design space of the continuous design space is  thus separated by the discrete candidate target set. One needs to decide the targets to be visited, and the order of targets. 

Although strict tree enumeration is a practical scheme for small-scale
space trajectory planning, in many such problems, exhaustive search is not tractable. Complexity of the search scheme quickly becomes intractable for the problems with even moderately sized candidate set. The  most challenging aspect of combinatorial optimization of space trajectory trajectory design are from the combinatorial part, i.e. the policy of selection of the candidate objects to construct the sequence.

To tackle the issues due to the problem's combintatoirial nature, one can discretize or branch the continuous design space into a set of separated continuous subspace. Each subspace are associated to a set of candidates using specific decision making policy. Graph based search strategies such as branch-and-bound tree searches are then used to search optimal sequence. \textcolor{red}{	It recursively splits the search space into smaller spaces associated to subset of the candidates, and performs a top-down recursive search through the tree of instances formed by the branch operation. Exploration and exploitation of the tree can be improved by introducing a knowledge based pruning criterion preventing nodes to be further expanded.} A notable example of tree searches with heuristic-based pruning, implementing problem knowledge, is the software STOUR of  automated design of trajectories with multiple fly-bys designed by Jet Propulsion Laboratory \cite{longuski1991automated}, which has been used in several important mission design works \cite{heaton2002automated}\cite{petropoulos2000trajectories}.

Studies in biological social action have always been a great source of inspiration of the mixed integer combinatorial space mission design. Attempts using bio-inspired  methodology, such as Ant Colony Optimization (ACO) \cite{ceriotti2010mga}, Genetic Algorithm (GA) \cite{Deb2007}\cite{Gad2011}\cite{izzo2014}, tree search strategies\cite{Petropoulos2014}  are made to advance automatic planning of the combinatorial search paradigms. Problem knowledge are used to define heuristics and tune of the internal parameters for different data sets or domains.

The ant colony meta heuristic and its variants have been successfully applied in a variety of  space optimziation problems. \textcolor{red}{The  Physarum and  the ACO variants model the discrete decision making problems into a decision graph, and unvisited node is stochasitcally selected  with probability} \cite{stuart2016design}\cite{simoes2017multi}\cite{ceriotti2010mga}.
Stuart et al. use a modified ACO suitable for smooth solution spaces to generate asteroid tours \cite{stuart2016design}. The route-finding ant colony optimization algorithm is incorporated into the automated tour-generation procedure for the sun–Jupiter Trojan asteroids. ACO routing algorithm and multi-agent coordination via auctions are  incorporated into a debris mitigation tour scheme design \cite{stuart2016application}. In  \cite{ceriotti2010mga}, an ant colony systems is used to solve MGAP design.  In \cite{simoes2017multi}, hybridization between Beam Search and the population-based ACO is proposed to tackle the combinatorial problem of finding the ideal sequence of bodies.

Approaches employing seeding based genetic algorithms are also proven effective in searching the optimal sequence  \cite{murakami2010approach}\cite{federici2019impulsive}.
Abdelkhalik et al. solve the multiple-flyby problem with impulsive chemical thrust  using a single genetic algorithm (GA) \cite{Abdelkhalik2012}. Nyew et al. proposed a structured-chromosome evolutionary algorithms for variable-size autonomous interplanetary trajectory planning optimization \cite{Nyew2015}. In \cite{Gad2011}, the problem is approached by using Hidden Genes Genetic Algorithms (HGGA) where each sequence is represented by a binary chromosome. \textcolor{red}{In multi-objective TSP with NSGA-II \cite{kdeb2007}, a representation scheme is used  to adapt the mixed integer decision variables to the NSGA-II framework. A 11-bit strings is designed to represent which and how many planets are used in the trajectory optimization. }

Other applications of metaheuristic algorithms in the path planning include: multi-objective particle swarm \cite{daneshjou2017mission}\cite{yu2015optimal}, simulated annealing \cite{cerf2015multiple}.  \cite{vasile2015incremental}\cite{di2017automatic} uses inspiration from the behavior of a simple amoeboid organism, the Physarum Polycephalum, that is endowed by nature with simple heuristics that can solve complex discrete decision making problems .

Some relatively recent works on discrete combinatorial problems of space trajectory design employ artificial intelligence and machine learning.  Tsirogiannis, George A proposed an  automatic planing approach using artificial intelligence and  machine learning \cite{tsirogiannis2012graph}. In \cite{hennes2015interplanetary},a heuristic-free approach using Monte Carlo Tree Search (MCTS) is proposed for the automated encounter sequence planning. Problem-specif modifications are made to adapt traditional MCTS to trajectory planning of the Rosetta and Cassini-Huygens interplanetary mission design. Four steps, namely: selection, expansion, simulation, and back-propagation, are performed iteratively till the stopping criteria are met \cite{browne2012a}.

It can be seen that, though numerous methods and algorithms are developed to tackle the space trajectory combinatorial optimization problem, either the graph based approach or the bio-inspired  and artificial intelligence based approaches,  the discretization  based searching strategy is used in natural.  Feasible space of the continuous parameters are discretized to construct branches and offspring associated to set of candidates. Policies such as bio-inspired social actions or population based evolutionary strategies  are then used to generate and select the offspring or branches to be explored. The procedure is iteratively implemented till the optimal sequence is obtained.  

\textcolor{red}{Some recent research works developed continuous representation techniques in combinatorial optimization. Pichugina et.al proposed a continuous presentation and functional extensions in combinatorial optimization. Functional representations of the Boolean set, general permutation, and poly permutation sets are derived\cite{Pichugina2016}. In \cite{2017Continuous}, a continuous representation technique in combinatorial optimization is proposed based on algebraic-topological features of the sets and properties of functions over them. } 

In this work, a new strategy for the sequence optimization is developed and tested. Different to the commonly used combinatorial optimization algorithms, a continuous gradient based optimizer for the space trajectory TSP  is developed. Candidates of the objects are associated to a set of continuous characteristic parameters, and a continuous mapping technique is developed to map the discrete or mixed integer combinatorial problem into a continuous quadratic-like optimization problem. Unlike the common combinatorial optimization algorithm, at each node of the sequence, a single candidate is selected and associated to a metric value, no branch and pruning operation is implemented. \textcolor{red}{A gradient based continuous optimizer is used to search optimal sequence.} Effectiveness of proposed method are demonstrated through a static TSP benchmark and multiple debris rendezvous problem .

The remainder of the paper is structured as follows. \textcolor{red}{ Section \ref{sec:problem} outlines the combinatorial space trajectory design and TSP of multiple debris rendezvous.} The continuous mapping approach is modelled and presented in Section \ref{sec:mapping}. \textcolor{red}{Section \ref{sec:test-static-TSP} presents tests of a static TSP benchmark using proposed method. Results of numerical simulation of the multiple debris rendezvous are presented in Section \ref{sec:results}}. Section \ref{sec:conclusion} concludes the paper and describe future work.

\section{Problem statement}\label{sec:problem}

\subsection{TSP in Space Mission Design}
The TSP in space mission design optimizes orbital transfer parameters and the sequence of objects to fly by. The solution space is defined by the launch date $T_0$ and times of flight (TOF) from one object to the subsequent one, $T_i$. The total cost of  the multiple transfer  can thus be described as

\begin{equation}\label{eq:cost function}
J =  \sum_{k=1}^N y(T_k,\mathcal{I}_k,\mathbf{d}_k)
\end{equation}
s.t.
\begin{equation}\label{eq:space craft}
\mathbf{s}_k = \mathbf{g}(\mathbf{s}_{k-1}) \\
\end{equation}
\begin{equation}\label{eq:obj-state}
\mathbf{x}_k = \mathbf{f}(\mathbf{x}_{k-1})\\
\end{equation}
\begin{equation}\label{eq:obj-observation}
\mathbf{z}_k = \mathbf{H}\mathbf{x}_k
\end{equation}
where $\mathcal{I}_k \in \mathbb{Z} $ represents the target body, $N$ is the number of targets. $\mathbf{d}_k $ is continuous design variables assoctaied to the $k$-th transfer,e.g.,relevant dates,  magnitude and orientation of impulse velocities  etc. $y(\cdot)$ is the cost function of the $k$-th transfer,  $\mathbf{s}$  are state of spacecraft e.g. position, velocity or orbital elements. $\textbf{x}_k$ and $\mathbf{z}_k$ are state and observation of the characteristic parameters of $\mathcal{I}_k$ respectively. $\mathbf{f}(\cdot)$, $\mathbf{g}(\cdot)$  are state equations of $\mathcal{I}_k$ and the spacecraft,  and $\mathbf{h}(\cdot)$ is observation equation of $\mathcal{I}_k$.  

The problem is a mixed integer combinatorial optimization problem. Optimization of the TSP is annoyed by its complex structured search spaces. One needs to explore not only the design space of the orbital transfer  parameters, but the discrete candidate set $\mathbf{I} \in \mathbb{Z}^N$.

The combination of the discrete set and continuous transfer parameters partitions the whole feasible solution space into a set of  separately distributed subspace.  Algorithms using different policies to generate and determine the optimal sequence are therefore developed. Though efforts are made to  improve the searching process, computational cost of the TSP solvers with even middle sized candidates is expensive. Since the time-dependent TSP is complicated by the discrete candidate set, if an equivalent continuous search space with continuous gradient for the TSP  can be  given, optimization of the problem will be much easier, and corresponding computational cost will be greatly reduced.

\subsection{Time-dependent TSP of multiple debris rendezvous }\label{sec:gtoc-9}

Consider a removal sequence problem of Sun-synchronous, Low-Earth-Orbit debris. Data set of the debris parameters are from the 9th Global Trajectory Optimisation Competition (GTOC9) \cite{izzo_dario_2018_1139022}. Inclinations and semi-major axes of the debris objects  are about 96$^0$ – 101$^0$ and  about 600 – 900 km larger than the Earth’s radius respectively. Eccentricity of the debris ranges from about 0.02 down to almost zero.  The task is to design a series of transfers to remove a set of orbiting debris.

\subsubsection{Problem statement}
Objective of the  multiple debris rendezvous is to  minimize 
\begin{equation}
J = \sum_{i=1}^N \Delta v_i  
\end{equation}
where $\Delta v_i$ is the cost of impulsive transfer to rendezvous the $i$-th debris.

After the rendezvous, the spacecraft stays in proximity of the debris for $5$ days for the removal operation and start a new transfer for the next rendezvous. Precession rates of the debris considering J2 coefficient given the orbital element $a,e,i,\Omega,\omega,M$ are computed as
\begin{align}\label{eqn:omega-rate}
\dot{\Omega} = -\frac{3}{2}J_2\left(\frac{r_{eq}}{p} \right)^2 n \cos i\\
\dot{\omega} = \frac{3}{4}J_2 \left( \frac{r_{eq}}{p} \right)^2 n (5\cos^2 i - 1) 
\end{align}
where $r_{eq}$ is mean radius of equator, $n = \sqrt{\frac{\mu}{a^3}}$ is mean motion, and $p$ is semilatus rectum and computed as $p=a(1-e^2)$.

Corresponding right ascension of the asecnding node, $\Omega$, argument of perigee, $\omega$ and the mean anomaly $M$  at epoch $t$ are
\begin{align}
\Omega - \Omega_0 = \dot{\Omega}(t - t_0) \\
\omega - \omega_0 = \dot{\omega}(t - t_0) \\
M - M_0 = n(t - t_0) 
\end{align} 
Initial values of the debris ephemerides, $[a_0,e_0,i_0,\Omega_0,\omega_0,M_0]$ and  initial epoch  $t_0$ are taken from the GTOC-9 database.

Cartesian form Ordinary Differential Equations (ODEs) of spacecraft's and debris' motion  is given by the following set of equations

\begin{align}
\ddot {\mathbf x} = - \frac{\mu x}{r^3} \left\{ 1 + \frac 32 J_2 \left( \frac {r_{eq}}r \right)^2 \left( 1 - 5 \frac {z^2}{r^2} \right) \right\}\\
\ddot {\mathbf y} = - \frac{\mu y}{r^3} \left\{ 1 + \frac 32 J_2 \left( \frac {r_{eq}}r \right)^2 \left( 1-  5 \frac {z^2}{r^2} \right) \right\}\\
\ddot {\mathbf z} = - \frac{\mu z}{r^3} \left\{ 1 + \frac 32 J_2 \left( \frac {r_{eq}}r \right)^2 \left( 3 - 5 \frac {z^2}{r^2} \right) \right\}
\end{align}
where  $r$ is the radius, $[x,y,z]$ and $[\dot{x},\dot{y},\dot{z}]$ are position and velocity respectively. Table \ref{Tab:constant} summarizes constants of the numerical integrator and computation of ephemerides of the debris.
\begin{table}[!h]
	\centering
	\caption{Values of problem constants }
	\begin{tabular}{lll}
		\hline\noalign{\smallskip}
		& Value         & Unit            \\
		\noalign{\smallskip}\hline\noalign{\smallskip}
		$\mu$       & 398600.4418   & km$^3$/sec$^2$  \\
		$J_2$       & 1.08262668e-3 &                 \\
		$r_{eq}$      & 6378.137      & km             \\          
		\hline\noalign{\smallskip}
	\end{tabular}\label{Tab:constant}
\end{table}
\subsubsection{Cost of near-circular impulsive orbital transfer}
Suppose the initial element of spacecraft at time $t_0$ is $\bm\sigma_0 = [a_0,e_0,i_0,\Omega_f,\omega_0,M_0]$, the $\Delta v$ for spacecraft fly to and rendezvous the target with orbital element $\bm\sigma_f = [a_f,e_f,i_f,\Omega_f,\omega_f,M_f]$ after time of flight $\Delta t$ can be approximately computed as

\begin{equation}\label{eq:dva}
\Delta V_{a} = \frac{1}{2}\frac{\Delta a}{a_0} V_0\\
\end{equation}
\begin{equation}\label{eq:dve}
\Delta V_{e} = \frac{1}{2}\Delta e V_0\\
\end{equation}
\begin{equation}\label{eq:dvi}
\Delta V_{i} = 2 V_0\sin\left(\frac{\Delta i}{2} \right)\\
\end{equation}
\begin{equation}\label{eq:dvo}
\Delta V_{\Omega} = \sin i_0 \Delta\Omega V_0 
\end{equation}
where $V_0$ is the spacecraft's velocity
\begin{equation}\label{key}
V_0 = \sqrt{\frac{\mu}{a_0} }
\end{equation}
and
\begin{align}
\Delta a = |a_0 - a_f|\\  
\Delta e = |e_0 - e_f|\\
\Delta i = |i_0 - i_f|\\
\Delta \Omega = |(\Omega_f + \dot{\Omega}_f \Delta t) - (\Omega_0 + \dot{\Omega}_0 \Delta t)|
\end{align}
where the drift rate of ascending node considering J2 impacts, $\dot{\Omega}_f$ and $\dot{\Omega}_0$ is computed using eq.\ref{eqn:omega-rate}.

The total cost of the transfer is 
\begin{equation}
\Delta V = \sqrt{\Delta V_a + \Delta V_e + \Delta v_i} + \Delta V_{\Omega}  
\end{equation}

\subsubsection{Reformulate the multiple rendezvous problem into a continuous optimization problem}

The optimization problem minimizes the total cost of $\Delta V_k $ and determine the optimal sequence $\mathcal{I}_k$. The mixed integer design vector of the optimization is $\{T_k,\mathcal{I}_k\}$, where $T_k$ is the ToF of the $k$-th transfer. To associate the integer $\mathcal{I}_k$ to continuous parameters, an expected debris $\textcolor{red}{\bar{\mathcal{I}}_k}$ is defined using   $\bm\mu_{k-1} = [\mu_a,\mu_e,\mu_i,\mu_o]_{k-1}$ and $\mathbf{\Sigma}_{k-1} = \text{diag}[\sigma_a,\sigma_e,\sigma_i,\sigma_o]_{k-1}$.  Initial value of $\bm\mu_{k-1}$ is set to differences of the orbital elements of $\mathcal{I}_k^*$ with respect to those of spacecraft at the $k-1$-th node. Observation of $\mathcal{I}_k$ is defined as $ \mathbf{z}_k = [\Omega_k]$, as the ascending node $\Omega_k$ is coupled closely to variation of the orbital parameters $a$,$e$, and $i$.

Design parameters of the continuous problem consist of   $\{T_k,\bm\mu_{k-1},\mathbf{\Sigma}_{k-1}\,\kappa_k\}$, where $\kappa_k$ is slack variable for the inequality constraint. The target $\mathcal{I}_k$ is selected using $\mathbf{z}_k$,  and  priori of $\bm\mu_{z,k}$ and $\mathbf{\Sigma}_{z,k}$ of $\textcolor{red}{\bar{\mathcal{I}}_k}$. Compute  $\Delta v_k$ to the debris $\mathcal{I}_k$, likelihood of the transfer $\mathcal{L}(\Delta v_k|\mathcal{I}_k)$ can then be computed with respect to $\mu_{v,k}$ and $\sigma_{v,k}$ of $\mathcal{I}_k^*$. Objective function value of the total cost with accumulated posterior constraints is  computed using eq.\ref{eqn:obj-slack-1}.

\section{Mapping discrete design parameters of TSP into continuous design space}\label{sec:mapping}

In state equation and observation equation of flyby target (eq.\ref{eq:obj-state}-eq.\ref{eq:obj-observation}), the state $\mathbf{x}_k$ and observation $\mathbf{z}_k$ is associated and varied with $\mathcal{I}_k$. Conventional TSP methods select the candidate $\mathcal{I}_k$ directly from the integer type data set  $\mathbf{I}$ using a specific policy.  Position and velocity of the object is then computed and substituent into the cost function to compute the total cost of the transfer. Selection and decision policy of the candidate $\mathcal{I}_k$ is analyzed and updated according to the total cost obtained.  The strategy works successfully in many space TSP optimization scenarios. However, as analyzed in preceding sections, discrete tree and string like structures have to be constructed and associated to the flyby targets. One have to examine each feasible subspace separately, no gradient of the trajectory variation of the spacecraft and flyby target is used to accelerate the searching process.

Instead of using integer type ID to represent the flyby target, in this work, a new set of continuous characteristic parameters are used to define the flyby target. Expected flyby target $\textcolor{red}{\bar{\mathcal{I}}_k}$ of the sequence are introduced, and a continuous probabilistic  metric with respect to expected $\textcolor{red}{\bar{\mathcal{I}}_k}$ is used to select the flyby target.  Accumulated probabilistic metric value of the  transfers are computed and set as constraints of the optimization.

\subsection{Associate integer ID to continuous characteristic parameters}

Suppose at the $k-1$-th node,  the state of spacecraft is $\mathbf{s}_{k-1}$. Initial value of the state of expected  $\textcolor{red}{\bar{\mathcal{I}}_k}$ is  $\mathbf{x}_{k-1}$ and normally distributed with $ \mathcal{N}(\bm\mu_{k-1},\mathbf{\Sigma}_{k-1})$. Expected value and covariance of the state of $\textcolor{red}{\bar{\mathcal{I}}_k}$ at the $k$-th node can be predicted and expressed as $\bm\mu_{k|k-1}$ and $\mathbf{\Sigma}_{k|k-1}$ respectively.

Given the predicted $\bm\mu_{k|k-1}$ and $\mathbf{\Sigma}_{k|k-1}$ , expected value of the observation of $\textcolor{red}{\bar{\mathcal{I}}_k}$ at the $k$-th node can be predicted as $\bm\mu_{z,k}$ and $\mathbf{\Sigma}_{z,k}$. A flyby target $\mathcal{I}_k$ can then be selected with respect to $\mathcal{I}_k^*$ using the characteristic parameter. Priori of the candidate $\mathcal{I}_k$, or the probability of $\mathcal{I}_k$ to be selected is given as
\begin{equation}\label{eq:priori}
P(\mathcal{I}_k|\mathcal{I}_k^*)   = \mathcal{N}(\mathbf{z}_k|\bm\mu_{z,k}, \mathbf{\Sigma}_{z,k}) 
\end{equation}   
where  $\mathbf{z}_k = h(\mathbf{x}_k)$ is characteristic parameters of $\mathcal{I}_k$ at the $k$-th node respectively.

Compute the probability level of each candidate $\mathcal{I}_k \in \mathbf{I}$  using eq.\ref{eq:priori}, and  save it in the data set of priori list for selection. Selection of $I_k$ is now associated to the probability level computed using a set of continuous parameters.

Response of the selection of $\mathcal{I}_k$ is also supposed to be Gaussian. Given the target $\mathcal{I}_{k}$, likelihood of the cost to  $\mathcal{I}_k$ can be compute as 
\begin{equation}\label{eq:likelihood}
\mathcal{L}(y_k |\mathcal{I}_k) = \mathcal{N}(y_k|\mu_{y,k}, \sigma_{y,k})
\end{equation}
where $y_k$ is the cost to $\mathcal{I}_k$, $\mu_y$ and $\sigma_y$ are expected value and variance of the cost of the $k$-th transfer.

With new introduced expected $\mathcal{I}_k^*$, selection of the target $\mathcal{I}_k$ and corresponding impacts on the cost of transfer can be rewritten  as a continuous problem.  Probability level of the transfer given the auxiliary characteristic parameters can be computed and measured. The original deterministic objective function can then be reformulated into an optimization problem that minimizes  the total cost of transfer while at the same time maximizes its accumulated probability level. This can be done by the following strategy. 

\subsection{Besysian Inference}

Eq.\ref{eq:likelihood} shows likelihood of the cost  $y_k$  given the priori of $\mathcal{I}_k$. With the priori and likelihood,  posterior of the transfer can be computed using Beaysian theorem
\begin{equation}
\pi(y_k) \propto \mathcal{L}(y_k |\mathcal{I}_k) P(\mathcal{I}_k|\mathcal{I}_k^*)  
\end{equation}

Accumulated logarithm of the posterior  for the TSP scheming is  
\begin{equation}
\begin{split}
\sum_{k=1}^N \log \pi(y_k) &= -\sum_{k=1}^N (\mathbf{z}_{k} - \bm\mu_{z,k}) \mathbf{\Sigma}_{z,k}^{-1} (\mathbf{z}_{k} - \bm\mu_{z,k})^T \\
& - \sum_{k=1}^N (y_k - \mu_{y,k}){\mathbf{\Sigma}}_{y,k}^{-1}(y - \mu_{y,k})^T + (\cdot)
\end{split}
\end{equation}
where $(\cdot)$ represents the logarithm term of
\begin{equation}\label{key}
 -\frac{1}{2} \log{\left(2 \pi\right)}- \frac{1}{2} \ln(|\mathbf{\Sigma}_{y,k}|) - \frac{1}{2}\ln(|\mathbf{\Sigma}_{z,k}|)- \frac{n_z}{2} \ln(2\pi) 
\end{equation}
and can be removed during the optimization if the variance of $y$ and $\mathbf{z}$ is small. The accumulated posterior includes the most current observation of the cost $y_k$, minimizing the value will allow us to move the samples in the prior to regions of high likelihood.
 
\begin{figure}
	\centering
	\includegraphics[width=1.2\textwidth]{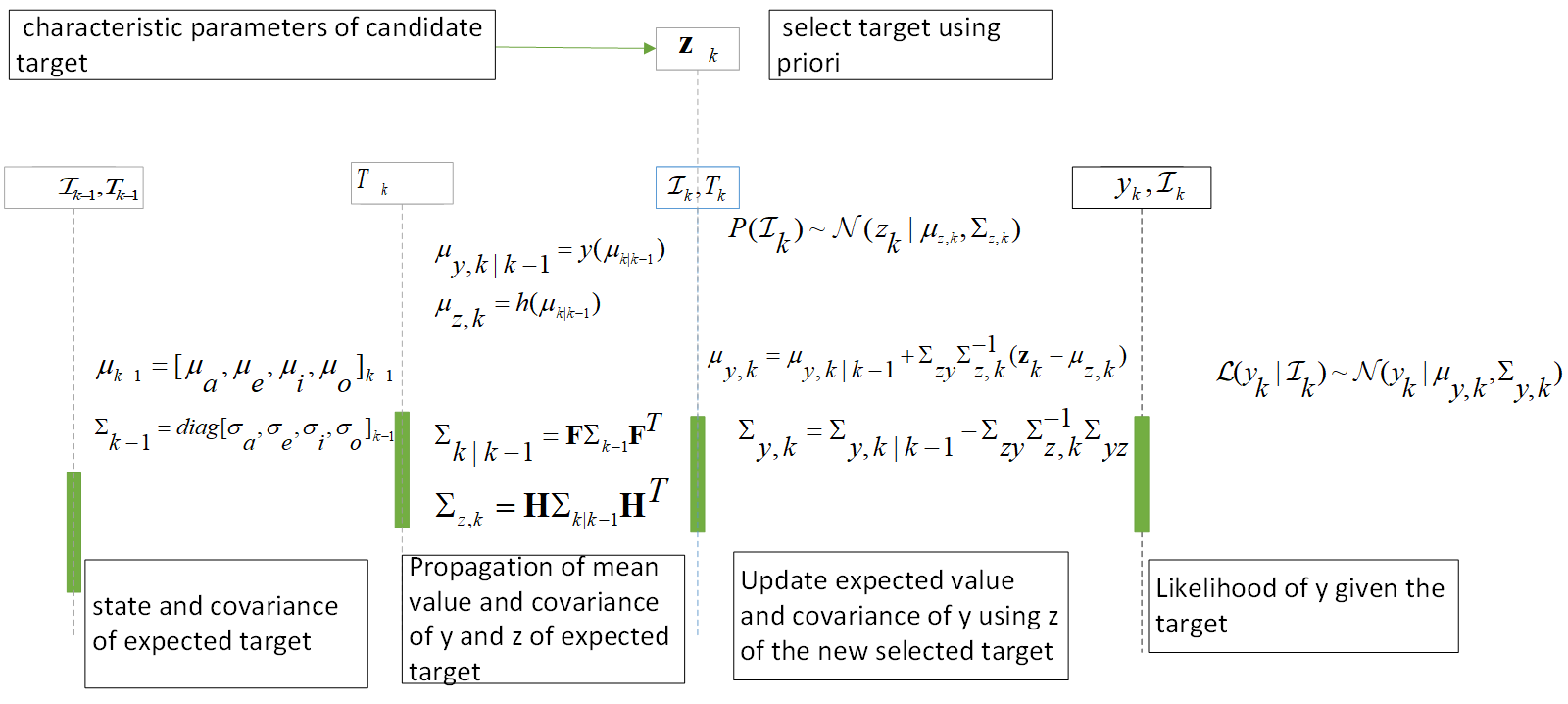}
	\caption{Priori and likelihood of selecting $\mathcal{I}_k$. $\mathcal{I}_k$ is selected using  characteristic parameters $\mathbf{z}$ with respect to $\mu_{\mathbf{z},k}$ and $\mathbf{\Sigma}_{z,k}$. Likelihood of $y$ is modeled and updated using the predicted mean value and covariance of $y$ given $\mathcal{I}_k$, and $\mathbf{z}_k$ of the selected $\mathcal{I}_k$.  }
	\label{fig:priori-likelihood}
\end{figure}

Figure \ref{fig:priori-likelihood} shows priori and likelihood of the expected $\mathcal{I}_k^*$ and flyby $\mathcal{I}_k$. First, an expected $\mathcal{I}^*_k$ is initialized and propagated to $t_k$. A flyby target $\mathcal{I}_k$ at $t_k$ is then selected from the data set using the priori. Corresponding cost to $\mathcal{I}_k$ is then computed. Posterior of the selection of $\mathcal{I}_k$ is computed using Bayesian theorem. Finally, a new objective function is defined to minimize the total cost of the transfer, while at the same to maximize the accumulated posterior.

\begin{equation}\label{eqn:new obj-1}
J = \min \sum_{k=1}^N y_k
\end{equation}
s.t.
\begin{equation}\label{eqn:new obj-2}
\begin{split}
\min  \Pi = &\sum_{k=1}^N (\mathbf{z}_{k} - \bm\mu_{z,k}) \mathbf{\Sigma}_{z,k}^{-1} (\mathbf{z}_k - \bm\mu_{z,k})^T \\
& + \sum_{k=1}^N (y_k - \bm\mu_{y,k})\mathbf{\Sigma}_{y,k}^{-1}(y_k - \bm\mu_{y,k})^T
\end{split}
\end{equation}

\subsection{Propagation of the expected values and covariance}

In the equations of priori and likelihood, several covariance of $\mathcal{I}_k$ and $\mathcal{I}_k^*$ , such as $ \mathbf{\Sigma}_{k|k-1}$, $\mathbf{\Sigma}_{z,k}$,$\mathbf{\Sigma}_{y,k}$ are involved and computed. The covariance are computed as follows.

Given initial value of $\mathcal{I}_k^*$ at the $k-1$th node, $\bm\mu_{k-1}$ and $\mathbf{\Sigma}_{k-1}$, expected value of $\mathcal{I}_K$ at the $k$-th node can be predicted as
\begin{equation}\label{eqn:propgation-1-mu}
\bm\mu_{k|k-1} = \mathbf{f}(\bm\mu_{k-1})
\end{equation}

Predict of the covariance of $\mathcal{I}_k^*$ at the $k$th node can be obtained using the first order approximation  
\begin{equation}\label{eqn:propgation-1-cov}
\mathbf{\Sigma}_{k|k-1} = \mathbf{F}_{k|k-1}\mathbf{\Sigma}_{k-1}\mathbf{F}_{k|k-1}^T
\end{equation}
where $\mathbf{F}_{k|k-1} $ is the Jacobian of state equation 
\begin{equation}
\mathbf{F}_{k|k-1} = \begin{bmatrix}
\frac{\partial f_1}{\partial x_1} & \frac{\partial f_1}{\partial x_2} & \cdots & \frac{\partial f_1}{\partial x_n}\\ 
\vdots &  & \ddots & \vdots \\ 
\frac{\partial f_n}{\partial x_1} & \frac{\partial f_n}{\partial g_2} & \cdots & \frac{\partial f_n}{\partial x_n}
\end{bmatrix}_{\mathbf{x} = \mathbf{x}_{k-1}}
\end{equation}

 Accordingly, $\bm\mu_{z,k|k-1}$ and $\mathbf{\Sigma}_{z,k|k-1}$ of expected $\mathcal{I}_k^*$ using the observation of characteristic parameters is
\begin{equation}\label{eqn:propgation-2-mu}
\bm\mu_{z,k} = \mathbf{h}(\bm\mu_{k|k-1})
\end{equation}
\begin{equation}\label{eqn:propgation-2-cov}
\mathbf{\Sigma}_\mathbf{z,k}= \mathbf{H}_{k} \mathbf{\Sigma}_{\mathbf{k|k-1}} \mathbf{H}_{k} ^T
\end{equation}

where $\mathbf{H}_{k|k-1} $ is the Jacobian of observation equation 

\begin{equation}
\mathbf{H}_{k} = \begin{bmatrix}
\frac{\partial h_1}{\partial x_1} & \frac{\partial h_1}{\partial x_2} & \cdots & \frac{\partial h_1}{\partial x_n}\\ 
\vdots &  & \ddots & \vdots \\ 
\frac{\partial h_n}{\partial x_1} & \frac{\partial h_n}{\partial x_2} & \cdots & \frac{\partial h_n}{\partial x_n}
\end{bmatrix}_{\mathbf{x} = \bm\mu_{k|k-1}}
\end{equation}

With the value of $\bm\mu_{k|k-1}$,$\mathbf{\Sigma}_{k|k-1}$, estimates of $\bm\mu_{y,k|k-1}$  $\mathbf{\Sigma}_{y,k|k-1}$ fly to expected flyby target $\mathcal{I}_k^*$ can be computed as 
\begin{equation}\label{eqn:propgation-4-mu}
\bm\mu_{y,k|k-1} = y(\bm\mu_{k|k-1})
\end{equation}
\begin{equation}\label{eqn:propgation-4-cov}
\mathbf{\Sigma}_{\mathbf{y},k|k-1}= \mathbf{Y}_{k} \mathbf{\Sigma}_{\mathbf{k-1}} \mathbf{Y}_{k} ^T
\end{equation}
where
\begin{equation}
\mathbf{Y}_{k|k-1} = 
\begin{bmatrix}
\frac{\partial y}{\partial x_1} & \frac{\partial y}{\partial x_2} & \cdots & \frac{\partial y}{\partial x_n}
\end{bmatrix}_{\mathbf{x}=\mathbf{x}_{k-1}} 
\end{equation}

Since the cost $y_k$ to actual $\mathcal{I}_k$ is computed and selected using the new observation $\mathbf{z}_k$, the following equation can be used to update covariance and expected value of  $y_k$ given  $\bm\mu_{z,k}$ and $\mathbf{\Sigma}_{z,k}$. 

Given jointly Gaussian $\mathbf{z}$,$\mathbf{y}$ , one have
 \begin{equation}
 \begin{bmatrix}
 y_k \\
 \mathbf{z}_k \\
 \end{bmatrix}
\sim \mathcal{N} \left ( 
\begin{bmatrix}
\bm\mu_{y,k}      \\
\bm\mu_{\mathbf{z,k}} 
\end{bmatrix},
 \begin{bmatrix}
\mathbf{\Sigma}_y                & \mathbf{\Sigma}_{y\mathbf{z}}      \\
\mathbf{\Sigma}_{y\mathbf{z}}    &\mathbf{\Sigma}_{\mathbf{z}}     
\end{bmatrix}
\right)  
\end{equation}
where
\begin{align}
\mathbf{\Sigma}_{z,k}  =  \mathbf{H}_k \mathbf{\Sigma}_{k|k-1}\mathbf{ H}_k^T\\
\mathbf{\Sigma}_{yy,k}  = \mathbf{Y}_k \mathbf{\Sigma}_{k|k-1}\mathbf{F}_k^T \\
\mathbf{\Sigma}_{yz,k}  = \mathbf{Y}_k \mathbf{\Sigma}_{k|k-1}\mathbf{H}_k^T \\
\mathbf{\Sigma}_{zy,k}  = \mathbf{H}_k \mathbf{\Sigma}_{k|k-1}\mathbf{F}_k^T 
\end{align}
Estimated expected value and covariance of $y_k$ given the observation of $\mathbf{z}_k$ can be computed as
\begin{equation}
\bm\mu_{y,k} =\bm\mu_{y,k|k-1} + \mathbf{\Sigma}_{yz}\mathbf{\Sigma}_{zz}^{-1}(\mathbf{z}_k - \bm\mu_{z,k|k-1})
\end{equation}
\begin{equation}
\mathbf{\Sigma}_{y,k} = \mathbf{\Sigma}_{y,k|k-1} - \mathbf{\Sigma}_{\mathbf{yz}}\mathbf{\Sigma}_{\mathbf{zz}}^-1\mathbf{\Sigma}_{\mathbf{z}y}
\end{equation}

\subsection{An alternative method for the covariance computation}

An alternative to the method of covariance update is using the scaled unscented transform (SUT) \cite{julier2002scaled}, which may give more exact results for large levels of the uncertainty but typically requires more computation time. With the method, predict of the mean vaule and covariance of $\mathbf{x}$ and $\mathbf{y}$ are computed as weighted sum of the sample function values of specified sigma points 
\begin{align}
\bar{\mathbf{x}} = \sum_{i=1}^r w_i \bm\chi_i \\
\mathbf{\Sigma}_{\mathbf{x}} = \sum_{i=1}^r w_i \left( \bm\chi_i - \bar{\mathbf{x}} \right) \left( \bm\chi_i - \bar{\mathbf{x}} \right) ^T \\
\bar{\mathbf{y}} = \sum_{i=1}^r w_i \mathbf{y}_i \\
\mathbf{\Sigma}_y = \sum_{i=1}^r w_i \left( \mathbf{y}_i - \bar{\mathbf{y}} \right) \left( \mathbf{y}_i - \bar{\mathbf{y}} \right) ^T \\
\mathbf{\Sigma}_{\mathbf{xy}} = \sum_{i=1}^r w_i \left( \bm\chi_i - \bar{\mathbf{x}} \right) \left( \mathbf{y}_i - \bar{\mathbf{y}} \right) ^T
\end{align}
where $w_i$ is the weight for the sampled points, $\bm\chi_i$ and $\mathbf{y}_i$ are sampled sigma points and corresponding function values. 

\subsection{Chi-square threshold}

In the searching procedure, the stochastic based computation and analysis  are used and implemented. Though the new objective function, eq.\ref{eqn:new obj-1}-eq.\ref{eqn:new obj-2} looks deterministic, the guess of $\mathcal{I}_k$ at each node is generated scholastically in essence (eq.\ref{eq:priori}). 

Consider the quadratic form cost function, eq.\ref{eqn:new obj-2}, since the guess of $\mathcal{I}_k$ is generated normally, the quadratic part is also randomly distributed because of the priori of $\mathcal{I}_k$,  and follow a chi-square distribution  
\begin{equation}
\left(\textbf{z}_k - \bm\mu_{z,k}\right) \mathbf{\Sigma}_{z,k} ^{-1} \left(\textbf{z}_k - \bm\mu_{z,k}\right) ^T + (y_k - \mu_{y,k})\mathbf{\Sigma}_{y,k}^{-1}(y_k - \mu_{y,k})^T \sim \chi^2_n  
\end{equation}
where $n = n_\mathbf{z} + n_y$ is degrees of freedom of chi-square distribution, and $n_\mathbf{z}$ and $n_y$ are dimensionality of the characteristic parameters and cost function value respectively.

\textcolor{red}{A chi-square constraint can then be given as}
\begin{equation}\label{eqn:obj-constr-2}
\begin{split}
\left(\textbf{z}_k - \bm\mu_{z,k}\right) \mathbf{\Sigma}_{z,k} ^{-1} \left(\textbf{z}_k - \bm\mu_{z,k}\right) ^T + (y_k - \mu_{y,k})\mathbf{\Sigma}_{y,k}^{-1}(y_k - \mu_{y,k})^T > \chi_n^2(\alpha)   
\end{split}
\end{equation}
where $\alpha$ is the confidence level to be met, e.g. 0.98.

Rewrite the constrained optimization problem into unconstrained form, one have
\begin{equation}\label{eqn:obj-slack-1}
\min_{\mathbf{D}} J = \sum_{k=1}^N \left[ y_k +  \kappa_i \max(0, c_k) \right] 
\end{equation}
\begin{equation}
\begin{split}\label{eqn:obj-slack-2}
c_k &= \left(\textbf{z}_k - \bm\mu_{z,k}\right) (\mathbf{\Sigma}_{z,k} + \mathbf{A}) ^{-1} \left(\textbf{z}_k - \bm\mu_{z,k}\right) ^T \\
&+ (y_k - \mu_{y,k})\sigma_{y,k}^{-1}(y_k - \mu_{y,k})^T
- \chi^2_n(\alpha) 
\end{split}
\end{equation}
where  $\mathbf{A} = \mathrm{diag}(\bm\epsilon)$, $\kappa_i$ is  slack variable, and $\mathbf{D}$ is design parameters of the new objective function 
\begin{equation}
\mathbf{D} = \{T_k,\mathbf{d}_k,\mu_{k-1}\,\mathbf{\Sigma}_{k-1},\kappa_k\} 
\end{equation}
in which $\mu_k$,$\mathbf{\Sigma}_k$ are the continuous auxiliary design variables of $\textcolor{red}{\bar{\mathcal{I}}_k}$, $\mathbf{d}_k$ is original continuous design variables of the TSP. 

The new  objective function incorporates statistical features of the expected $\textcolor{red}{\bar{\mathcal{I}_k}}$, and the original cost function, and is reformulated into a continuous optimization problem associated to a probability space. The objective is to search the optimal solution of sequence that has lowest cost at maximum probability level. A gradient based optimizer can  be used to search iteratively the characteristic parameters of $\textcolor{red}{\bar{\mathcal{I}}_k}$ till  the chi-square threshold is achieved.

\subsection{Searching optimal TSP sequence using gradient based optimizer}

In the selection using priori,  $\mathcal{I}_k$ is randomly selected with respect to the spacecraft and expected $\textcolor{red}{\bar{\mathcal{I}}_k}$. Though the operation is similar to many probabilistic heuristic algorithms for the combinatorial optimization, several features of the proposed algorithm distinguish it from them and make it high efficient in searching the optimal sequence. 

With the new strategy, continuous gradient information of the accumulate logarithm probabilistic level of $\mathcal{I}_k$ with respect to  expected $\textcolor{red}{\bar{\mathcal{I}}_k}$ can be computed and used. No branch operation is required during the search. In the population based methods in which the probabilistic heuristics  with the string like structured samples are used and analyzed to to obtain heuristics for the new individuals. In the new algorithm, only one candidate is selected at each node. The computation is implemented sequentially.  No population based sampling and seeding  is required.

The algorithm works as follows.First,  initial values of the characteristic parameters, $\bm\mu_{k-1}$ and $\mathbf{\Sigma}_{k-1}$ are given to initialize the expected $\textcolor{red}{\bar{\mathcal{I}}_k}$.  State equation of the spacecraft and flyby target are used to propagate the characteristic parameters. Select $\mathcal{I}_k$ from the candidate set using the observations with the prirori of  $\bm\mu_{z,k}$ and $\mathbf{\Sigma}_{z,k}$. Compute expected value and covariance of the cost, $\bm\mu_{y,k}$ and $\mathbf{\Sigma}_{y,k}$. Slacked inequality constrained objective function value is computed using eq.\ref{eqn:obj-slack-1}. Finally, a gradient based sequential quadratic programming (SQP) solver is  used to search iteratively the optimal continuous design parameters $\{T_k\}$,$\{\mathbf{d}_k\}$, $\{\bm\mu_{k-1}\}$, $\{\mathbf{\Sigma}_{k-1}\}$ and the slack variable $\{\kappa_k\}$. Algorithm I summarizes the searching procedure.

\begin{algorithm}
	\caption{Optimization of Time-dependent TSP }
	\label{alg:tsp}
	\begin{algorithmic}[1]
		
		\Procedure{OPTIMIZETSP}{$\mathbf{T}$,$\mathbf{I}$}
    \State \textbf{Input:}  Initial values of $\{T_k\}$,$\{\mathbf{d}_k\}$, $\{\bm\mu_{k-1}\}$, $\{\mathbf{\Sigma}_{k-1}\}$, $\{\kappa_k\}$ and $n_{iter}$  
	\State \textbf{Output:}  Optimal $\mathcal{I}_k^*$,  $T_k^*$, and $\mathbf{d}_k^*$   
		
		\State $ iter = 0$   
		\While{$iter < n_{iter}$ or $\nabla J > \epsilon$}  \Comment{ the optimal solution?}
		\State $iter = iter + 1$
		\For {$k$th node}
		\State $\mathbf{x}_k$ $\leftarrow$ $\mathbf{f}(\mathbf{x}_{k-1})$
		\State  propagate $\bm\mu_{k|k-1}$,$\mathbf{\Sigma}_{k|k-1}$ of $\bar{\mathcal{I}}_k$  \Comment{eq.\ref{eqn:propgation-1-mu}- eq.\ref{eqn:propgation-1-cov}} 
		\State Propagate $\bm\mu_{z,k}$,$\mathbf{\Sigma}_{z,k}$ of $\textcolor{red}{\bar{\mathcal{I}}_k}$   \Comment{eq.\ref{eqn:propgation-2-mu} -eq.\ref{eqn:propgation-2-cov}}
		\State Select an object $\mathcal{I}_k$ from  $\mathbf{I}$  using priori with respect to $\textcolor{red}{\bar{\mathcal{I}}_k}$ \Comment{ eq.\ref{eq:priori}}
		\State Characteristic parameters $\mathbf{z}_k$ of $\textcolor{red}{\bar{\mathcal{I}}_k}$
		\State Parameter of new cost function $\mu_{y,k},\sigma_{y,k}$ of $\textcolor{red}{\bar{\mathcal{I}}_k}$     \Comment{eq.\ref{eqn:propgation-4-mu}-eq.\ref{eqn:propgation-4-cov}}
		\State  $y_k$ $\leftarrow$ transfer cost of $\mathcal{I}_k$ 
		\EndFor
		\State $J$ $\leftarrow$ objective function with the slack variable  \Comment{eq.\ref{eqn:obj-slack-1} - eq.\ref{eqn:obj-slack-2}}
		\State Update $\{T_k\}$,$ \{\mathbf{\mathbf{d}}_k\} $, $\{\bm\mu_{k-1}\}$, $\{\mathbf{\Sigma}_{k-1}\}$, $\{\kappa_k\}$ using gradient based optimizer
		\EndWhile  
		\State Return $\{T_k\}, \{\mathcal{I}_k\}$ and $\{\mathbf{d}_k\}$ 
		\EndProcedure
		
	\end{algorithmic}
\end{algorithm}

\section{Test of Static TSP}\label{sec:test-static-TSP}
\textcolor{red}{
Consider a test case of static TSP. The benchmark TSP consists of 14 points. Table \ref{tab:TSP-Benchmark-Position} shows coordinates of the 14 points \cite{aravindseshadri2021}. The objective of the problem is to find an optimal travel route such that each points is visited one and only once with the least possible distance traveled.  Aravind Seshadri resolved this problem using a stochastic combinatorial optimization algorithm, Simulated Annealing (SA) algorithm \cite{aravindseshadri2021}.}

\begin{table}[H]
	\label{tab:TSP-Benchmark-Position} 
	\centering
	\caption{Coordinates of the points of TSP benchmark}
	\begin{tabular}{lll}
		\hline
		ID    & $x$ & $y$  \\
		\hline
1  & 16.470 &	96.100\\
2  & 16.470 &	94.440\\
3  & 20.090 &	92.540\\
4  & 22.390 &	93.370\\
5  & 25.230 &	97.240\\
6  & 22.000 &	96.050\\
7  & 20.470 &	97.020\\
8  & 17.200 &	96.290\\
9  & 16.300 &	97.380\\
10 & 14.050 &	98.120\\
11 & 16.530 &	97.380\\
12 & 21.520 &	95.590\\
13 & 19.410 &	97.130\\
14 & 20.090 &	94.550\\
		\hline
	\end{tabular}
\end{table}

\textcolor{red}{Different to the problem of optimizing tour of rendezvousing the debris, in test of benchmark, design variables are expected value and variance of $x$ and $y$ to be travelled to the next point. Candidate points are selected via analysing  probabilities of possible candidates with the data set of parameter $\mu_x$,$\mu_y$,$\sigma_x$,$\sigma_y$. Table \ref{Tab:Optimal-TSP-benchmark} shows the optimal solution of the benchmark using SA algorithm and corresponding design variables of $\mu_x$,$\mu_y$. Both starting point and final point are set to point 13. The minimum total cost of TSP is 30.8785. }

\begin{table}[H]
	\centering
	\label{Tab:Optimal-TSP-benchmark}
	\caption{Optimal solutions of the benchmark }
	\begin{tabular}{crr} 
		\hline
ID  &     $\mu_x$   &   $\mu_y$ \\
		\hline
13  &  			    &           \\
7  &  			1.06&        -0.11  \\
12  &  			1.05&        -1.43  \\
6  &  			0.48&         0.46  \\
5  &  			3.23&         1.19  \\
4  &  		   -2.84&        -3.87  \\
3  &  			2.3&         -0.83  \\
14  &               0&         2.01  \\
2  &           -3.62&        -0.11  \\
1  &               0&         1.66  \\
10  &           -2.42&         2.02  \\
9  &            2.25&        -0.74  \\
11  &            0.23&            0  \\
8  &            0.67&        -1.09   \\
13 &                &    \\
		\hline
	\end{tabular}
\end{table}

\textcolor{red}{Unlike the time-dependent TSP whose covariance of the node is predicted by numerical propagating differential equations, in the static TSP, variance of the position $i$, $\sigma_i$ given prior variance of position $i$, $\sigma_i^-$ (the variance parameter should be optimized) and variance at point $i-1$, $\sigma_{i-1}$, variance at point $i$ can be computed as 
\begin{equation}
\sigma_{i}^2 = (\sigma_i^-)^2 + \sigma_{i-1}^2 + 2\rho(\sigma_i^-)\sigma_{i-1}  
\end{equation}
where $\rho\in [0,1]$ is the correlation factor. Candidate point $i$ is selected using the distribution parameter $[\mu_{x,i}, \sigma_{x,i}]$ and$  [\mu_{y,i}, \sigma_{y,i}] $.}

\textcolor{red}{Corresponding variance of distance $r_i = \sqrt{x_i^2 + y_i^2}$ of $i$th point can be computed as 
\begin{equation}
\sigma_{r,i}^2 = \frac{(\sigma_{x,i}^2\mu_{x,i}^2 + \sigma_{y,i}^2\mu_{y,i}^2)}{(\mu_{x,i}^2 + \mu_{y,i}^2)}; % succeed
\end{equation}
Therefore, the dsign variables for each point consists of the distribution parameters  $\mu_x,\mu_y,\sigma_x,\sigma_y$, correlation factor $\rho$, and slack variable $\kappa$ for the constraint to be met. Table \ref{Tab:Design-Space-Benchmark} summarizes design variables of one point on the travel route. Total number of the design variables is therefore $7\times 13 = 91$. }

\begin{table}[H]
	\centering
	\label{Tab:Design-Space-Benchmark}
	\caption{Design variables of the static TSP }
	\begin{tabular}{llll}
		\hline
		Parameters  & Initial value & Lower bound & Upper bound  \\
		\hline
		$\mu_x$     &      -         &  -8.0       &  8.0        \\
		$\mu_y$     &      -         &   8.0       &  6.5        \\
		$\sigma_x$  &      4.0       &   0.1       &  6.0        \\
		$\sigma_y$  &      4.0       &   0.1       &  6.0        \\
		$\rho_x$    &      0.2       &   0.0       &  1.0        \\
		$\rho_y$    &      0.2       &   0.0       &  1.0        \\
		$\kappa$    &  50            & 0.01        &  300        \\
		\hline
	\end{tabular}
\end{table}

In this test, we provide two different constraints of the route. One is based on Chi-square threshold as the preceding section does. The objective function is set to minimize the total cost with chi-square constraints at each point
\begin{equation}\label{key}
\min J = \sum_{i=1}^{14} r_i + \kappa_i \max(0, c_i)
\end{equation}
where $c_i$ is the inequality constraints at each point
\begin{equation}
\begin{split}
c_i &= \left(\textbf{z}_i - \bm\mu_{z,i}\right) \mathbf{\Sigma}_{z,i}  ^{-1} \left(\textbf{z}_i - \bm\mu_{z,i}\right) ^T \\
&+ (r_i - \mu_{r,i})\sigma_{r,i}^{-1}(r_i - \mu_{r,i})^T
- \chi^2_n(\alpha) 
\end{split}
\end{equation} 
and $\textbf{z}_i = [x_i,y_i], \bm\mu_{z,i}= [\mu_{x,i},\mu_{y,i}]$, and $\mathbf\Sigma_{z,i} = \text{diag}([\sigma_{x,i},\sigma_{y,i}])$ respectively.

Another one is based on  Maximum A Posteriori (MAP) and can be written as
\begin{multline}
\min J = \sum_{i=1}^{14} r_i + \log(|\mathbf{\Sigma}_{z,i}|) + \log \sigma_{r,i}+\\
 (\mathbf{z}_i - \bm\mu_{z,i})\mathbf\Sigma_{z,i}^{-1}(\mathbf{z}_i - \bm\mu_{z,i})^T + (r_i - \bm\mu_{r,i})\sigma_{r,i}^{-1}(r_i - \mu_{r,i})^T  
\end{multline}
In the equation of logarithm of posterior, items that have no influence to the optimal solution are removed.

Table \ref{Tab:Solution-ABC-of-Benchmark} shows two groups of results obtained with different parameter settings. Initial values of $\mu_{x,i}$ and $\mu_{y,i}$ are given deterministically or randomly within a specific range near the optimal design variables, where $\mu_{x,i}^*$ and $\mu_{y,i}^*$ are the optimal design variables of $\mu_{x,i}$ and $\mu_{y,i}$,  $\text{rnd}\in[0,1]$ is a uniformly distributed random number in the interval (0,1).  \footnote{Source code of the tesst of benchmark can be available at   https://github.com/Liqiang-Hou/Travelling-Salesman-Problem-of-Space-Trajectory-Design}.

\begin{table}[H]
	\centering
	\label{Tab:Solution-ABC-of-Benchmark}
	\caption{Solutions with different initial parameter settings }
	\begin{tabular}{ccccl} 
		\hline
		Solution &Constraint &Initial value of $[\mu_{x,i},\mu_{y,i}]$& Route & Cost  \\ 
		\hline
		A        & MAP &  $[\mu_{x,i}^*+1.0 , \mu_{y,i}^*-0.8]$                          &    [13,7,12,6,5,4,3,14,2,1,10,9,11,8,13]    & 30.8785      \\
		B        & Chi-square   & $[-2 + 4\text{rnd} ,-2 + 4\text{rnd}]$                         &  [13,7,12,6,5,4,3,14,2,1,8,11,9,10,13]      &  31.567     \\
		\hline
	\end{tabular}
\end{table}

It can be seen that if the initial values of design variables is within a certain  neighbourhood  of the optimal ones, one can obtain the global optimal solution directly using the single constraint of MAP. Regarding solution B, though the initial values are away from the optimal initial design parameters and randomly selected,  the solutions obtained are still quite close to the global optimal solution.  

With the numerical results, one can see that though the solution obtained is not necessarily global optimal for its gradient-based searching mechanism, however, because of the Chi-square inequality constraint at each point, convexity of the problem is improved (the inequality constraints can be written into equally into quadratic form equality constraints with Lagrange multiplier), and the more constraints the problem has, the more convex the problem can be, and the easier the optimizer could search the optimal route.

The proposed optimizer for the mixed integer TSP is designed assuming that the nodes are Gaussian correlated. The assumption looks reasonable in the test cases presented in the manuscript, but maybe too strict for some other TSP problems. We note that in the tests of TSP benchmark, within a specific range, given uniformly random distributed initial parameters, one can still find the near optimal solution. Therefore, if a more sophisticated distance measure other than  Euclidean distance, and some more sophisticated correlation model can be given, generality of proposed approach could be improved. New sampling techniques and continuous stochastic constraints thresholds for the algorithm might help improve the algorithm’s generality too. Such techniques for improving generality of proposed method in the future work are under considered.

\section{Numerical simulation}\label{sec:results}

Consider now a time-dependent multiple debris rendezvous problem. Ephemerids of the debris are from GTOC-9 website \cite{izzo2014}. Table \ref{tab:debris-statistic} summarizes  distribution characteristics of the  debris.

\begin{table}[ht]
	\caption{Debris data set}
	\label{tab:debris-statistic} 
	\begin{tabular}{llllll}
\hline\noalign{\smallskip}
		\multicolumn{1}{c}{}    & \multicolumn{1}{c}{unit} & \multicolumn{1}{c}{mean value} & standard deviation & maximum  & minimum    \\ 
\noalign{\smallskip}\hline\noalign{\smallskip}
		semi-major axis, $a$    & km                       & 7131.6                         & 58.503             & 7274.0   & 6996.1     \\
		eccentricity, $e$       &                          & 0.0071438                      & 0.004901           & 0.019318 & 0.00013122 \\
		inclination, $inc$      & deg                      & 98.415                         & 0.84304            & 101.07   & 96.236     \\
		Ascending node $\Omega$ & deg                      & 172.93                         & 103.88             & 347.76   & 7.6191     \\ 
\hline\noalign{\smallskip}
	\end{tabular}
\end{table}

\subsection{Optimal solution of JPL}

Table \ref{tab:optim-jpl} shows the optimal solution of JPL using an improved genetic based evolutionary algorithm (GA). Initial branch and bound, beam search, and seeding combinatorial algorithms (Genetic Algorithm, GA)are used to search the optimal chain of rendezvous. A new structure genome, and heuristic adjustment are used to improve the GA. The optimal solution to remove the total 123 debris consist of 10 sequences.  Table \ref{tab:optim-jpl} shows one  of the optimal sequences, its optimal transfer duration and costs of the sequence \cite{Petropoulos2014} .

	\begin{table}[ht]
\caption{Optimal solution of JPL}
\label{tab:optim-jpl} 
\begin{tabular}{ll}
\hline
	Parameter               & Value                                                                                                                     \\
\hline
	Start MJD2000           & 23557.18                                                                                                                  \\
	End MJD2000             & 23821.03                                                                                                                  \\
	Number of  objectives   & 14                                                                                                                        \\
	Debris ID               & 23,55,79,113,25,20,27,117,121,50,95,102,38,97                                                                             \\
	Transfer Duration, days & \begin{tabular}[c]{@{}l@{}}24.86,24.98,22.42,24.99,0.29,10.63,25.00,\\  2.70,1.51,1.41,24.67,24.31,5.86\end{tabular}      \\
	$\Delta v$, m/s         & \begin{tabular}[c]{@{}l@{}}161.8,139.2,65.8,208.2,115.2,\\ 300.1,564.9,78.3,105.0,233.3,\\ 453.5,340.4,300.8\end{tabular}\\
\hline
\end{tabular}
	\end{table}

	\subsection{Optimization of the sequence using continuous gradient optimizer}

Consider the  multiple debris rendezvous problem. Initial epoch and the first debris to remove of the mission is set to	$\mathrm{MJD}_0 = 23557$ and  $\mathcal{I}_0 = 23$ as those of JPL.

Similar to JPL's solution, debris are selected and computed from the database of 123 debris, removal operation of the debris after each transfer is set to 5 days. As analyzed in preceding section, the expected value and covariance of expected debris $\mathcal{I}_i^*$, and ToF of each transfer are set to the parameters to be optimized. Table \ref{tab:lb-ub-test} shows the initial settings of the design parameters of each transfer. For each transfer, there are 10 design variables need to be optimized. The total number of the design parameters is therefore $ 14\times 10 =140$. The lower bound and upper bound of $\mu_{\Omega}$ are set to $-8^0$ and $8^0$ respectively as the maximum value of $\Delta \Omega$ obtained by varying $\Delta a$, $\Delta i$ and $\Delta e$ of the debris  is about $5^0$, corresponding upper bound of $\sigma_{\Omega}$ is set to $8^0$ for the same reason.

\begin{table}[!htb]
\caption{Parameter settings of the design parameters}
\label{tab:lb-ub-test} 
		\begin{tabular}{lllll}
\hline
			Parameters        & unit & lower bound & upper bound & initial value \\ \hline
			T                 & days & 0.5         & 25          & 20            \\
			$\mu_{a}$         & km   & -150        & 150         & 0             \\
			$\mu_{e}$         &      & -1.0e-3     & 1.0e-3      & 0             \\
			$\mu_{i}$       & deg  & -1.5        & -1.5        & 0             \\
			$\mu_{\Omega}$    & deg  & -8.0        & 8.0         & 0             \\
			$\sigma_{a}$      & km   & 5.0         & 50          & 30            \\
			$\sigma_{e}$      &      & 1.0e-4      & 1.0e-3      & 5.0e-4        \\
			$\sigma_{i}$    & deg  & 0.1         & 1.0         & 0.5           \\
			$\sigma_{\Omega}$ & deg  & 0.1         & 8           & 5             \\ 
			$\kappa$ 		  &      & 1.0e-3      & 300         & 50             \\ 
\hline
		\end{tabular}
	\end{table}

The finite difference method is used to compute the Jocabian for covariance propagation. The Matlab$\textsuperscript{\textregistered}$ sequential quadratic programming (SQP) solver, fmincon is used to search the optimal sequence and ToF of the transfer.  Table \ref{tab:optimzer-settings} shows parameter settings of the finite difference, fmincon optimizer and Chi-square threshold of the optimization.

\begin{table}[!htb]
	\caption{Parameter settings of the optimization}
	\label{tab:optimzer-settings} 
	\begin{tabular}{lllll}
		\hline\noalign{\smallskip}
		Parameter                                                   &              & unit & value       \\
		\hline\noalign{\smallskip}
		\multirow{4}{*}{Finite difference for covariance propagation} & $\delta a$      & km   & 1.0e-2      \\
		& $\delta e$      &      & 1.0e-4      \\
		& $\delta i$      & deg  & 1.0e-3      \\
		& $\delta  \Omega$ & deg  & 1.0e-3      \\
		Chi-square threshold                                       &              &      & $\chi^2_2(0.98)$ \\
		\multirow{2}{*}{Parameter settings of fmincon optimizer}                           & TolX         &      & 1.0e-20     \\
		& TolFun       &      & 1.0e-20   \\ 
		\hline\noalign{\smallskip}
	\end{tabular}
\end{table}
\subsection{Numerical results and analysis}
\begin{figure}[ht]
	\centering
	\includegraphics[width=0.8\textwidth]{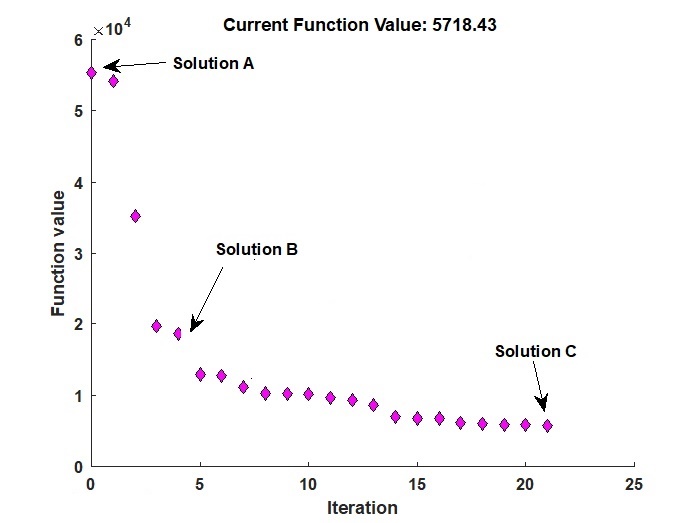}
	\caption{Convergence history }
	\label{fig:solution-convergence}
\end{figure}

First, a sequence of transfers with fixed ToF is optimized. The value of ToF is set to 20 days for each transfer between the debris. Figure \ref{fig:solution-convergence} shows convergence history  and three solutions of proposed algorithm. Table \ref{tab:optim-tsp-1-seq}-Table \ref{tab:optim-tsp-3-var} shows corresponding sequence, initial mean value and variance of the three solutions. Penalty of the cost function and velocity of each transfer are listed in the tables too.  With the algorithm, it takes couples of iterations to reduce the penalties of transfer cost sufficiently close to zero. Variance and mean value of the  sequence orbital elements converge also to a sufficient low level.    

The gradient based optimizer, fmincon of Matlab$\textsuperscript{\textregistered}$ is again used  to refine the ToF of each transfer of the sequence. Table \ref{tab:optim-tsp-4} summarizes the optimal solution of the sequence. A sequence of total $\Delta V = 3337.0$m/s si obtained, while the total cost of JPL is 2966.5m/s. In this test,   analytical approximation formulas are used to compute the transfer cost,  while in JPL's results, a global multiple shooting based local optimization with numerical orbital propagator is used. Considering accuracy level of the two models used, the costs of two test sequences can be seen at similar level.  

From the tables, it can be seen that with proposed algorithm, both the mean value and variances of the orbital parameters converge gradually to  appropriate values that associate to the IDs of targets. However, the variance of semi-major axis, $\sigma_a$ remains almost unvaried. This is due to the impact of the semi-major axis to the $\Delta V$. Given a variation of semi-major axis $\Delta a = 30$km,  the $\Delta V$ varied due to  $\Delta a$  is about  15m/s, while for the ascending node and inclination, this value could be be up to more than 100m/s if $\Delta \Omega$ and $\Delta i$ are set to $1^0$ respectively (eq.\ref{eq:dva}-eq.\ref{eq:dvo}). Therefore, the value of  $\sigma_a$ varies significantly less than other design parameters during the optimization. 

One can also see that, though the design parameter of ToF and sequence can be optimized directly following the algorithm, in this test, a two-step strategy is used to search the optimal solution of the sequence. First, an optimal sequence with fixed-ToF is determined, a local optimizer is then used to search optimal ToF and $\Delta V$ for each transfer. Since the ToFs are not considered in the first step, some actual optimal solutions may be missed during the optimization. The reason that the two-step strategy is used is related to the covariance and cost function propagation strategy in the optimization.   In this test, the first-order approximation for the covariance propagation is used.  The state equation of orbital elements are highly nonlinear. Transition matrix of the states of orbital elements can not be given explicitly, the finite difference is therefore used to compute indirectly the covariance due to varied states. For the similar reason, it is hard to propagate the covariance by integrating the Ricatti equation in which an explicit state transition matrix is required. Therefore, the fixed ToF  optimization is used to avoid the complex varied time-step covariance propagation. With the approximation, one needs to tune  carefully the finite difference of the covariance propagation, and optimality of the solutions is dependent on the finite difference used.  

We note the recent sigma-point technique for covariance propagation. With sigma-point technique, covariance of the state and observation can be propagated directly using the nonlinear state and observation equations.  Such a work that directly optimizes the sequence using derivative-free covariance propagation  will be included in future work. A state equation with position and velocity will also be considered, in which the state transition matrix can be explicitly listed.  In this test, part of the orbital elements related design variables, e.g. $\sigma_{\Omega}$ is quite insensitive to the varied cost function values. If design variables of position vectors are introduced, sensitivity of the design variables to the varied cost function will be improved.
   
\begin{table}[ht]
\caption{Solution A: the sequence}
\label{tab:optim-tsp-1-seq} 
	\begin{tabular}{lllll}
\hline\noalign{\smallskip}
Penalty & Tof(days) & Delta V(km/s) & $\mathcal{I}_{k-1}$ & $\mathcal{I}_k$ \\
\hline\noalign{\smallskip}
		3833.6 & 20 & 0.33356 & 23  & 55  \\
		3630.9 & 20 & 0.33285 & 55  & 113 \\
		4109   & 20 & 0.19693 & 113 & 121 \\
		2515.4 & 20 & 0.59675 & 121 & 117 \\
		3921   & 20 & 0.30597 & 117 & 79  \\
		3694   & 20 & 0.34758 & 79  & 20  \\
		3720.5 & 20 & 0.3394  & 20  & 25  \\
		4107.9 & 20 & 0.18808 & 25  & 27  \\
		3848.2 & 20 & 0.26756 & 27  & 84  \\
		4202.1 & 20 & 0.71163 & 84  & 83  \\
		4038.2 & 20 & 0.96796 & 83  & 50  \\
		3948.4 & 20 & 0.6506  & 50  & 118 \\
		3282.8 & 20 & 0.41598 & 118 & 95  \\
		2364.2 & 20 & 0.63115 & 95  & 87  \\
\hline\noalign{\smallskip}
	\end{tabular}
\end{table}

\begin{table}[ht]
\caption{Solution A: mean value and variance}
\label{tab:optim-tsp-1-var} 
	\begin{tabular}{llllllll}
\hline\noalign{\smallskip}
		\begin{tabular}[c]{@{}l@{}}$\mu_a$\\ (km)\end{tabular} &
\begin{tabular}[c]{@{}l@{}}$\mu_i$\\ (deg)\end{tabular} &
\begin{tabular}[c]{@{}l@{}}$\mu_o$\\ (deg)\end{tabular} &
\begin{tabular}[c]{@{}l@{}}$\sigma_a$\\ (km)\end{tabular} &
\begin{tabular}[c]{@{}l@{}}$\sigma_i$\\ (deg)\end{tabular} &
\begin{tabular}[c]{@{}l@{}}$\sigma_o$\\ (deg)\end{tabular} &
\begin{tabular}[c]{@{}l@{}}$\mu_v$\\ (km/s)\end{tabular} &
\begin{tabular}[c]{@{}l@{}}$\sigma_v$\\ (km/s)\end{tabular} \\
\hline\noalign{\smallskip}
		0 & 0 & 0 & 30 & 1 & 3 & 0 & 0.33976 \\
		0 & 0 & 0 & 30 & 1 & 3 & 0 & 0.34934 \\
		0 & 0 & 0 & 30 & 1 & 3 & 0 & 0.34984 \\
		0 & 0 & 0 & 30 & 1 & 3 & 0 & 0.35219 \\
		0 & 0 & 0 & 30 & 1 & 3 & 0 & 0.3603  \\
		0 & 0 & 0 & 30 & 1 & 3 & 0 & 0.35758 \\
		0 & 0 & 0 & 30 & 1 & 3 & 0 & 0.34931 \\
		0 & 0 & 0 & 30 & 1 & 3 & 0 & 0.36707 \\
		0 & 0 & 0 & 30 & 1 & 3 & 0 & 0.33538 \\
		0 & 0 & 0 & 30 & 1 & 3 & 0 & 0.3365  \\
		0 & 0 & 0 & 30 & 1 & 3 & 0 & 0.35832 \\
		0 & 0 & 0 & 30 & 1 & 3 & 0 & 0.32785 \\
		0 & 0 & 0 & 30 & 1 & 3 & 0 & 0.34322 \\
		0 & 0 & 0 & 30 & 1 & 3 & 0 & 0.34184 \\
\hline\noalign{\smallskip}
	\end{tabular}
\end{table}

\begin{table}[!htb]
\caption{Solution B: the sequence}
\label{tab:optim-tsp-2-seq} 
	\begin{tabular}{llllll}
\hline\noalign{\smallskip}
Penalty & Tof(days) & Delta V(km/s) & $\mathcal{I}_{k-1}$ & $\mathcal{I}_k$ \\
\hline\noalign{\smallskip}
		193.76  & 19.984 & 0.33333  & 23  & 55   \\
		526.91  & 19.984 & 0.33294  & 55  & 113  \\
		1.9254  & 19.984 & 0.096311 & 113 & 79   \\
		193.8   & 19.984 & 0.29176  & 79  & 121  \\
		741.25  & 19.983 & 0.42998  & 121 & 117  \\
		632.04  & 19.984 & 0.49426  & 117 & 57   \\
		611.02  & 19.99  & 0.69273  & 57  & 20   \\
		737.28  & 20.003 & 0.47951  & 20  & 27   \\
		217.48  & 20.003 & 0.26803  & 27  & 84   \\
		0.96473 & 20.002 & 0.11276  & 84  & 83   \\
		30.334  & 20.002 & 0.265    & 83  & 50   \\
		94.584  & 20.003 & 0.25041  & 50  & 118  \\
		412.72  & 20.005 & 0.34214  & 118 & 25   \\
		910.09  & 20.006 & 0.54754  & 25  & 95   \\
\hline\noalign{\smallskip}
	\end{tabular}
\end{table}

\begin{table}[!htb]
\caption{Solution C: mean value and variance}
\label{tab:optim-tsp-2-var} 
	\begin{tabular}{llllllll}
\hline\noalign{\smallskip}
		\begin{tabular}[c]{@{}l@{}}$\mu_a$\\ (km)\end{tabular} &
\begin{tabular}[c]{@{}l@{}}$\mu_i$\\ (deg)\end{tabular} &
\begin{tabular}[c]{@{}l@{}}$\mu_o$\\ (deg)\end{tabular} &
\begin{tabular}[c]{@{}l@{}}$\sigma_a$\\ (km)\end{tabular} &
\begin{tabular}[c]{@{}l@{}}$\sigma_i$\\ (deg)\end{tabular} &
\begin{tabular}[c]{@{}l@{}}$\sigma_o$\\ (deg)\end{tabular} &
\begin{tabular}[c]{@{}l@{}}$\mu_v$\\ (km/s)\end{tabular} &
\begin{tabular}[c]{@{}l@{}}$\sigma_v$\\ (km/s)\end{tabular} \\
\hline\noalign{\smallskip}
		80.241 & 0.20271  & 0.25162  & 29.691 & 0.59927 & 2.5142 & 0.057949 & 0.1238  \\
		68.961 & 0.14367  & 0.1659   & 29.777 & 0.33223 & 2.4502 & 0.052913 & 0.10902 \\
		6.9494 & 0.019553 & 0.071687 & 29.997 & 0.4815  & 2.6018 & 0.011746 & 0.14774 \\
		56.518 & 0.13509  & 0.16085  & 29.701 & 0.65134 & 2.4081 & 0.03901  & 0.12259 \\
		117.79 & 0.30928  & 0.37205  & 29.353 & 0.45844 & 2.4847 & 0.084661 & 0.11743 \\
		171.14 & 0.39852  & 0.45827  & 29.379 & 0.52076 & 2.4762 & 0.15569  & 0.11906 \\
		219.99 & 0.54225  & 0.63295  & 29.624 & 0.30007 & 2.46   & 0.1493   & 0.11256 \\
		143.58 & 0.37154  & 0.43753  & 29.338 & 0.49681 & 2.4938 & 0.10586  & 0.12098 \\
		52.673 & 0.13962  & 0.17868  & 29.706 & 1.2963  & 2.5114 & 0.036761 & 0.18234 \\
		10.65  & 0.029357 & 0.09848  & 29.985 & 0.48589 & 2.3233 & 0.013828 & 0.12152 \\
		56.707 & 0.14237  & 0.16824  & 29.805 & 0.67806 & 2.5316 & 0.052932 & 0.13629 \\
		46.089 & 0.07849  & 0.0956   & 29.919 & 0.37909 & 2.5339 & 0.043194 & 0.11532 \\
		76.403 & 0.15493  & 0.17917  & 29.827 & 0.36931 & 2.4881 & 0.061738 & 0.11272 \\
		179.85 & 0.44911  & 0.55832  & 29.189 & 0.49114 & 2.518  & 0.12415  & 0.11769 \\
\hline\noalign{\smallskip}
	\end{tabular}
\end{table}
% Please add the following required packages to your document preamble:
% \usepackage{booktabs}
% Please add the following required packages to your document preamble:
% \usepackage{booktabs}

\begin{table}[!htb]
\caption{Solution C: sequence}
\label{tab:optim-tsp-3-seq} 
	\begin{tabular}{llllllll}
		\hline\noalign{\smallskip}
Penalty & Tof(days) & Delta V(km/s) & $\mathcal{I}_{k-1}$ & $\mathcal{I}_k$ \\
		\hline\noalign{\smallskip}
		0       & 20 & 0.33365  & 23  & 55  \\
		0.43601 & 20 & 0.33282  & 55  & 113 \\
		0       & 20 & 0.096486 & 113 & 79  \\
		0.76409 & 20 & 0.29294  & 79  & 121 \\
		0.33184 & 20 & 0.42921  & 121 & 117 \\
		0       & 20 & 0.4898   & 117 & 57  \\
		0       & 20 & 0.68889  & 57  & 20  \\
		0.59349 & 20 & 0.47947  & 20  & 27  \\
		0.85352 & 20 & 0.26742  & 27  & 84  \\
		0       & 20 & 0.11133  & 84  & 83  \\
		0.94618 & 20 & 0.2685   & 83  & 50  \\
		0.90029 & 20 & 0.25063  & 50  & 118 \\
		0.64101 & 20 & 0.34123  & 118 & 25  \\
		0       & 20 & 0.54873  & 25  & 95  \\
		\hline\noalign{\smallskip}
	\end{tabular}
\end{table}

\begin{table}[!htb]
\caption{Solution C: mean value and variance}
\label{tab:optim-tsp-3-var} 
	\begin{tabular}{llllllll}
		\hline\noalign{\smallskip}
		\begin{tabular}[c]{@{}l@{}}$\mu_a$\\ (km)\end{tabular} &
		\begin{tabular}[c]{@{}l@{}}$\mu_i$\\ (deg)\end{tabular} &
		\begin{tabular}[c]{@{}l@{}}$\mu_o$\\ (deg)\end{tabular} &
		\begin{tabular}[c]{@{}l@{}}$\sigma_a$\\ (km)\end{tabular} &
		\begin{tabular}[c]{@{}l@{}}$\sigma_i$\\ (deg)\end{tabular} &
		\begin{tabular}[c]{@{}l@{}}$\sigma_o$\\ (deg)\end{tabular} &
		\begin{tabular}[c]{@{}l@{}}$\mu_v$\\ (km/s)\end{tabular} &
		\begin{tabular}[c]{@{}l@{}}$\sigma_v$\\ (km/s)\end{tabular} \\
		\hline\noalign{\smallskip}
53.706  & 0.16283   & 1.656    & 19.688 & 0.31928 & 0.36124 & 0.22875  & 0.016381 \\
41.364  & 0.19485   & 1.4922   & 19.268 & 0.56979 & 1.8124  & 0.23917  & 0.10725  \\
-12.849 & 0.28146   & 2.4948   & 18.348 & 0.36068 & 0.83147 & 0.46402  & 0.039063 \\
40.8    & -0.068711 & -0.52225 & 20.063 & 0.68693 & 2.1619  & 0.15941  & 0.13708  \\
77.781  & 0.16261   & 1.2189   & 18.847 & 0.56092 & 1.5637  & 0.15396  & 0.083502 \\
138.13  & 0.031337  & 1.3298   & 18.45  & 0.35088 & 0.92499 & 0.075472 & 0.015769 \\
146.61  & 0.17877   & 1.1054   & 19.041 & 0.47357 & 1.6726  & 0.12666  & 0.070588 \\
96.794  & 0.15679   & 1.0286   & 18.981 & 0.6832  & 2.1857  & 0.11309  & 0.14058  \\
32.133  & 0.041496  & 0.3307   & 19.331 & 0.75819 & 3.4516  & 0.031035 & 0.26005  \\
-30.017 & 0.4219    & 2.8362   & 13.339 & 0.51202 & 0.79282 & 0.5873   & 0.057609 \\
37.748  & 0.15001   & 0.88379  & 19.886 & 0.78621 & 7.9787  & 0.13354  & 1.1676   \\
25.887  & 0.10616   & 0.77631  & 19.589 & 0.82071 & 4.1815  & 0.11944  & 0.38334  \\
47.493  & 0.14112   & 1.0868   & 19.514 & 0.60612 & 2.1292  & 0.15945  & 0.12676  \\
120.86  & 0.11869   & 0.84214  & 18.596 & 0.5103  & 1.7112  & 0.065413 & 0.028208 \\
		\hline\noalign{\smallskip}
	\end{tabular}
\end{table}

\begin{table}[ht]
\caption{Optimal sequence}
\label{tab:optim-tsp-4} 
	\begin{tabular}{@{}lll@{}}
		\hline\noalign{\smallskip}
		Parameter & unit & Value               \\                                                                                     \hline\noalign{\smallskip}
		ID        &      & 23, 55 , 113, 79 , 121, 117, 57 , 20, 27, 84, 83, 50, 118, 25,95 \\
		Tof &
		days &
		\begin{tabular}[c]{@{}l@{}}8.5433 , 24.966 , 19.166 , 24.964 , 24.994 , 24.989 , 24.996 ,\\  23.834 , 21.521 , 5.0216 , 14.567 , 24.986 , 24.989\end{tabular} \\
		$\Delta V$ &	m/s &
		\begin{tabular}[c]{@{}l@{}}160.38 , 351.14 , 77.365 , 261.65 , 411.83 ,  245.97 , 370.18 , \\475.27 , 184.97 , 77.182 , 192.64 , 255.65 , 272.76\end{tabular} \\ 
		Total $\Delta V$ &	m/s & 3337.0 \\
		\hline\noalign{\smallskip}
	\end{tabular}
\end{table}

\section{Conclusions}\label{sec:conclusion}
The time-dependent TSP problem is hard to resolve due to its complicated design space structure. Conventional combinatorial algorithms such as search tree and seeding combinatorial algorithms use selection and regenerating operators based on the integer type variables. Efficiency of the searching procedure is highly dependent on the policy used for the integer based operation.

Different to the conventional TSP methods, the new algorithm resolves the problem using a gradient based search.  A probabilistic mapping is introduced to  associated the integer ID to continuous characteristic parameters. Parameter estimation techniques such as  covariance propagation are used to determine the optimal characteristic parameters. No tree or graph based operation is involved. Though the optimizer is designed using stochastic techniques, the searching procedure is implemented deterministically.\textcolor{red}{Chi-square threshold checking operation is introduced into the optimization, and used to improve convexity of the problem.}

Another contribution of the strategy is from the Baseysin based heuristics. Stochastic meta heuristics are commonly used in  prevalent TSP algorithms. In these algorithms, individuals and branches are generated and sampled using probabilistic heuristic, but few of them uses Bayesian based response analysis directly in the operations due to the discrete and mixed integer nature. With the continuous mapping and Bayesian  analysis, objective function of the sequence is reformulated into a continuous quadratic-like search problem. As far as the authors know, the research study is the first trial in space mission design to employ a continuous gradient optimizer to resolve mixed integer time-dependent space TSP.

The method is designed stochasticly in essence, a robust design of TSP considering uncertainty impact using the proposed method can be considered  in future work. Extensions and application of the proposed algorithm to the TSP with multiple objectives are well worth trying too. Applications of the continuous mapping and Byesian mechanism to the multi-target flybys and multi-gravitational assist are also considered.

\begin{acknowledgements}
This work is supported by National Natural Science Foundation of China, No.U20B2056,U20B2054.

\end{acknowledgements}

% Authors must disclose all relationships or interests that 
% could have direct or potential influence or impart bias on 
% the work: 
%
% \section*{Conflict of interest}
%
% The authors declare that they have no conflict of interest.

% BibTeX users please use one of
%\bibliographystyle{spbasic}      % basic style, author-year citations
\bibliographystyle{spmpsci}      % mathematics and physical sciences
%\bibliographystyle{spphys}       % APS-like style for physics
%\bibliography{}   % name your BibTeX data base

\bibliography{reference}

\end{document}